\def\n${\begin{equation}}
\def\N${\end{equation}}

\documentstyle{article}
\parskip=\baselineskip

\begin{document}

\begin{center}
{\Large Two-variable Wiman-Valiron theory and PDEs\\ Erratum}
\\\vspace{6mm}
P.C. Fenton
and
John Rossi
\end{center}

\setcounter{equation}{23}

\noindent

\noindent Abstract: {\small This note corrects Example 3.2 in a  published paper with the same title.\\
 Mathematics Subject Classification 32A15}

\noindent Typographical errors and a careless comment on the Cauchy-Kovalevskaya theorem make Example 3.2 in our published paper [1] somewhat confusing.  Below is a corrected version. 
\\[.2cm]
\noindent {\bf Example 3.2}

\noindent Consider the $n$-th order linear PDE
\begin{equation} \label{lin}
\sum_{i+j\leq n} P_{i,j}f^{i,j}=0,
\end{equation} where the $P_{i,j}$ are polynomials in two complex
variables. We assume that $f$ is an entire solution of (\ref{lin}). Clearly,
as in Example 3.1, every such solution is transcendental.

To the best of our knowledge there have been no order estimates of
entire solutions of (\ref{lin}). Our method can often obtain such
results. To simplify matters, let us take the second order equation
\begin{equation} \label{2nd}
f^{1,1}=Pf,
\end{equation}
where $P$ is a polynomial, and proceed as in Example 3.1. Using (15), (\ref{2nd}) becomes
$$
A {\zeta}_1 {\cal F}^{1,0}+ B {\zeta}_2 {\cal F}^{0,1}+ C
{\zeta}_{1}^2 {\cal F}^{2,0}+D {\zeta}_1{\zeta}_2 {\cal F}^{1,1}+
E{\zeta}_{2}^2{\cal F}^{0,2}={\zeta}_{1}^4{\zeta}_{2}^2
P({\zeta}_1^2{\zeta}_2,{\zeta}_1{\zeta}_2^2){\cal F}, $$ where $A$,
$B$, $C$, $D$ and $E$ are constants. Using (11), we
obtain \begin{equation} \label{N1N2} A{\cal N}_1 + B{\cal N}_2 + C
{\cal N}_{1}^2 +D {\cal N}_1{\cal N}_2 + E{\cal
N}_{2}^2=(1+o(1)){\zeta}_{1}^4{\zeta}_{2}^2P({\zeta}_1^2{\zeta}_2,{\zeta}_1{\zeta}_2^2),\end{equation}
and therefore \begin{equation} \label{newton} {{\cal N}^*}^2 \geq
(c+o(1))|{\zeta}_{1}|^4|{\zeta}_{2}|^2|P({\zeta}_1^2{\zeta}_2,{\zeta}_1{\zeta}_2^2)|,\end{equation}
where $c$ is a positive constant. As in Example 3.1, this implies
that\\
 $\rho({\cal F})\geq 3+3d/2$, where $d$ is the degree of $P$, and
thus $\rho(f)\geq 1+d/2$.

\noindent We remark also that the reference after equation (7) should read [2, p. 228].\\[.3cm]

\noindent
{\large References}\\

\noindent
1.  Fenton, P.C. and Rossi, John, Two Variable Wiman-Valiron Theory and PDEs, Ann. Acad. Sci. Fenn. Math 35 (2010), 571-580.\\[.3cm]

\noindent University of Otago, Dunedin, New Zealand\\ (e-mail:pfenton@maths.otago.ac.nz)\\[.2cm]
Virginia Tech, Blacksburg, VA, USA\\ (e-mail:rossij@math.vt.edu)

\end{document}